\documentclass{elsart3-1}

\usepackage{amssymb, amsmath}

\usepackage[english,francais]{babel}

\newtheorem{theorem}{Theorem}[section]
\newtheorem{lemma}[theorem]{Lemma}
\newtheorem{proposition}[theorem]{Proposition}
\newtheorem{corollary}[theorem]{Corollary}

\newtheorem{conjecture}[theorem]{Conjecture}

\renewcommand{\theequation}{\arabic{equation}}
\def\og{\leavevmode\raise.3ex\hbox{$\scriptscriptstyle\langle\!\langle$~}}
\def\fg{\leavevmode\raise.3ex\hbox{~$\!\scriptscriptstyle\,\rangle\!\rangle$}}
\def\EMdash{\leavevmode\hbox to 7.5mm{\vrule height .63ex depth -.59ex width 5.4mm\hfill}}

\let\epsilon=\varepsilon

\begin{document}

\begin{frontmatter}




%
\selectlanguage{francais}
\title{\textbf{La Structure de $A$-Module induite par un  $A$-Module de Drinfeld de Rang 2
sur un corps fini}}

\vspace{-2.6cm} \selectlanguage{english}
\title{\textbf{The $A$-Module Structure Induced by a Drinfeld $A$-Module of Rank 2
over a Finite Field}}



\author[a]{Mohamed-Saadbouh MOHAMED-AHMED}
\ead{mohamed-saadbouh.mohamed-ahmed@univ-lemans.fr}

\address[a]{D\'epartement de Math\'ematiques, Universit\'e du Maine, Avenue Olivier Messiaen,
72085 Le Mans Cedex 9, France}

\selectlanguage{francais}
\begin{abstract}
Soit $\Phi $ un $\mathbf{F}_{q}[T]$-module de Drinfeld de rang
$2$, sur un corps fini $L$, extension de degr\'e $n$ d'un corps
fini $\mathbf{F}_{q}$. Soit $P_{\Phi }(X)=$ $X^{2}-cX+\mu
P^{m}$ (o\'u $c \in \mathbf{F}_{q}[T],$ $\mu $ est un \'el\'ement non nul de $%
\mathbf{F}_{q}$, $m$ est le degr\'e  de l'extension $L$ sur  $%
\mathbf{F}_{q}[T]/P,$ et $P$ est la $\mathbf{F}_{q}[T]$-caract\'eristique de $%
L $ et $d $ le  degr\'e du  polyn\^ome $P$) le  polyn\^ome
caract\'eristique du Frobenius $F$ de $L$. On s'int\'eressera \'a
la structure de $\mathbf{F}_{q}[T]$-module fini  $L^{\Phi }$
induite par $\Phi $ sur $L$. Notre r\'esultat principal est le
parfait analogue du th\'eor\`eme de  Deuring ( voir
\cite{Deuring}) pour les courbes elliptiques
:\ \ \ soit $M=\frac{\mathbf{F}_{q}[T]}{I_{1}}\oplus \frac{\mathbf{F}_{q}[T]}{%
I_{2}}$, o\`u $I_{1}=(i_{1})$ et $\ \ I_{2}=(i_{2})$\ ( $i_{1}$,
$i_{2}$ sont deux  polyn\^omes  de $\mathbf{F}_{q}[T]$) tels que :
$i_{2}\mid (c-2)$. Il existe alors un $\mathbf{F}_{q}[T]$-module
de Drinfeld ordinaire $\Phi $ sur $L$ de rang $2$ tel que:
$L^{\Phi }$ $\simeq M$.

{\it Pour citer cet article: Mohamed-saadbouh.Mohamed-Ahmed, C.
R. Acad. Sci. Paris,
 Ser. I ... (...).}

\vskip 0.5\baselineskip \selectlanguage{english}
\noindent{\bf Abstract} \vskip 0.5\baselineskip \noindent Let
$\Phi $ be a Drinfeld $\mathbf{F}_{q}[T]$-module of rank $2$,
over a finite field $L$. Let $P_{\Phi }(X)=$ $X^{2}-cX+\mu P^{m}$
($c$
an element of $\mathbf{F}_{q}[T],$ $\mu $ be a non-vanishing element of $%
\mathbf{F}_{q}$, $m$  the degree of the extension $L$ over the field $%
\mathbf{F}_{q}[T]/P,$ and $P$  the $\mathbf{F}_{q}[T]$-characteristic of $%
L $ and $d $  the degree of the polynomial $P$) the
characteristic polynomial of the Frobenius $F$ of $L$. We will be
interested in the structure of finite $\mathbf{F}_{q}[T]$-module
$L^{\Phi }$ induced by $\Phi $ over $L$. Our main result is
analogue to that of Deuring ( see \cite{Deuring} ) for elliptic
curves
:\ \ \ Let $M=\frac{\mathbf{F}_{q}[T]}{I_{1}}\oplus \frac{\mathbf{F}_{q}[T]}{%
I_{2}}$, where $I_{1}=(i_{1})$,$\ \ I_{2}=(i_{2})$\ ( $i_{1}$,
$i_{2}$ being two polynomials of $\mathbf{F}_{q}[T]$) such that :
$i_{2}\mid (c-2)$. Then there exists an ordinary Drinfeld
$\mathbf{F}_{q}[T]$-module $\Phi $ over $L$ of rank $2$ such that
: $L^{\Phi }$ $\simeq M$. {\it To cite this article:
Mohamed-Saadbouh Mohamed-Ahmed , C. R. Acad. Sci. Paris, Ser. I
... (...).}
\end{abstract}

\end{frontmatter}

\section{Introduction}

 let $K$ a no empty global field of characteristic $p$ (
 namely a rational functions field of one indeterminate over a
 finite field ) together with a constant field, the finite field
 $\mathbf{F}_{q}$ with $p^{s}$ elements. We fix one place of $K$,
 denoted by $\infty ,$ and call $A$  the ring of regular elements
 away from the place $\infty $. Let $L$ be a commutator field of
 characteristic $p$,  $\gamma :A\rightarrow L$ be a ring
 $A$-homomorphism. The kernel of this $A$-homomorphism is denoted by $P.$
 We put $m$ =$[L,$ $A/P]$, the extension degree of $L$ over
 $A/P$, and $d= deg P$.

 We denote by $L\{\tau \}$  the
 polynomial ring of $\tau $, namely the Ore polynomial ring, where $\tau $ is the Frobenius of
 $\mathbf{F}_{q}$ with the usual addition and where
 the product is given by the commutation rule: for every $%
 \lambda \in L\mathit{,}$ we have $\tau \lambda =\lambda ^{q}\tau
 $. A Drinfeld $A$-module  $\Phi :A \rightarrow
 L\{\tau \}$ is a non trivial ring homomorphism and a non trivial embedding
 of $A $ into $L\{\tau \}$ different from $\gamma $.
 This homomorphism $\Phi $, once defined, define an $A$-module structure over the $A$-field $L$, noted \ $%
 L^{\Phi }$, where the name of a Drinfeld $A$-module for a
 homomorphism $\Phi $. This structure of $A$-module depends on
 $\Phi $ and, especially, on his rank, for more information see
    \cite{B.Angles},  \cite{Gosse}, and \cite{Drinfeld2}.


We will be interested in  a Drinfeld $A$-module structure
$L^{\Phi }$ in the case
of rank 2, and we will prove that for an ordinary Drinfeld $\mathbf{F}%
_{q}[T] $-module, this structure is always the sum of two cyclic and finite $%
\mathbf{F}_{q}[T]$-modules : $\frac{A}{I_{1}}\oplus \frac{A}{I_{2}}$ where $%
I_{1}=(i_{1})$ and $I_{2}=(i_{2})$ such that $i_{1}$ and $i_{2}$
are two ideals of $A,$ which verifies $i_{2}\mid i_{1}$. Let $P_{\Phi }(X)=$ $X^{2}-cX+\mu P^{m},$ such that $\mu \in \mathbf{F}%
   _{q}^{\ast } $, and $c$ $\in A$, the characteristic polynomial of $\Phi$.We will
show that $\chi
_{\Phi }=I_{1}I_{2}=(P_{\Phi }(1))$, so if we put \ $i=$pgcd$%
(i_{1},i_{2}), $ then : $i^{2}\mid P_{\Phi }(1)$. \ We will give
an analogue of Deuring theorem for elliptic curves :
\begin{theorem}
Let $M=\frac{A}{I_{1}}\oplus \frac{A}{I_{2}}$, where
$I_{1}=(i_{1})$ ,$\ \ I_{2}=(i_{2})$\ and such that : $\
i_{2}\mid i_{1}, $ $i_{2}\mid (c-2)$. Then there exists an
ordinary Drinfeld $\ \mathit{A}$-module $\Phi $ over $L$ of rank
$2$, such that: $L^{\Phi }\simeq M$.
\end{theorem}






\section{Structure de $A$-module de Drinfeld $L^{\Phi }$}
 The Drinfeld $A$-module of rank $2$ is of the form $\ \Phi (T)=a_{1}+a_{2}\tau +a_{3}\tau
 ^{2}$, where $\ a_{i}\in L$, $1\leq i\leq 2$ and $a_{3}\in L^{\ast}$. Let $\Phi $ and $\Psi $ be two Drinfeld modules over an
$A$-field \ $L$. A morphism from $\Phi $ to $\Psi $ over $L$ is an element $ p(\tau )\in L%
     \mathit{\{\tau \}}$ such that
     $p \Phi _{a}=\Psi_{a} p\quad \mbox{ for all }
     a\in A.$
   A non-zero morphism is called an isogeny. We note that this is
   possible only between two Drinfeld modules with the same rank.
   The set of all morphisms forms an $A-$module denoted by Hom$_{L}(\Phi ,\Psi ).$

    In particular, if $\Phi $ =$\Psi $ the $L$-endomorphism ring End$%
    _{L}\Phi =$Hom$_{L}(\Phi ,\Phi $) is a subring of \ $L\{\tau \}$ and an $A$%
    -module contained in $\Phi (A)$. Let $\overline{L}$ be a fix algebraic closure of
$L$, $\Phi _{a}(\overline{L}):=\Phi \lbrack a](\overline{L})=\{x\in \overline{L}%
,\Phi _{a}(x)=0\}$, and  $\Phi _{P}(\overline{L})=\cap _{a\in P}\Phi _{a}(%
                                                  \overline{L}).$
     We say that $\Phi $ is supersingular if and only if the $A$-module constituted by a $P$%
     -division points $\Phi _{P}(\overline{L})$ is trivial,
     otherwise  $\Phi $ is said an ordinary module, see
     \cite{Gosse}.

Let $\Phi $ be a Drinfeld $A$-module of rank $2,$ \ over a finite
field $L$
and let $P_{\Phi }$ his characteristic polynomial, $P_{\Phi }(X)=$ $X^{2}-cX+\mu P^{m},$ such that $\mu \in \mathbf{F}%
   _{q}^{\ast } $, and $c$ $\in A,$ where $\ \deg c\leq
   \frac{m.d}{2}$ by the Hasse-Weil analogue in this case.
    Let $\chi $ be the Euler-Poincar\'e characteristic ( i.e. it is an
   ideal from $A$). So we can speak about the ideal $\chi (L^{\Phi})$, denoted henceforth by $\chi _{\Phi }$, which is by definition
   a divisor of $A$, corresponding for the elliptic curves to a
   number of points of the variety over their basic
   field.
   About the $A$-module structure $%
L^{\Phi }$, we have the following result :

\begin{proposition}
The Drinfeld $A$-module $\Phi $ give a finite $A$-module
structure $L^{\Phi
} $, which is on the form $\frac{A}{I_{1}}\oplus \frac{A}{I_{2}}$ where $%
I_{1}$ and $I_{2}$ are two ideals of $A,$ such that: $\chi _{\Phi
}=I_{1}I_{2}.$
\end{proposition}

We put $I_{1}=(i_{1})$ and $\ \ I_{2}=(i_{2})$ ( $i_{1}$ and
$i_{2}$ two unitary polynomials in $A$).

Let  $i=$ pgcd $(i_{1},i_{2})$, it is clear by the Chinese lemma,
that the no cyclicity of the $A$-module $L^{\Phi }$, needs that \
$I_{1}$ and $I_{2}$ are not a prime between them, that means that
$i\neq 1,$ and since the relation $\chi _{\Phi }=I_{1}I_{2}$, we
will have : $\ i^{2}\mid P_{\Phi }(1) $ $\ (\chi _{\Phi
}=(P_{\Phi }(1))$).

In all the next of this paper, the condition above, will be
considered
verified, and more precisely we suppose that $I_{2}\mid I_{1}$ $\ ($i.e $%
: i_{2}\mid i_{1})$ otherwise $L^{\Phi }$ is a cyclic $A$-module
and can be writing on this form $A/\chi _{\Phi }.$

\begin{proposition}
If $L^{\Phi }\simeq \frac{A}{I_{1}}\oplus \frac{A}{I_{2}}$, then
$i_{2}\mid c-2$.
\end{proposition}

Proof : We know that the  $A$-module structure $L^{\Phi }$ is
stable by the endomorphisme  Frobenius $F$ of $L$. We choose a
basis for $ A /\chi _{\Phi },$ for which  the $A$-module
$L^{\Phi }$ will be generated by $(i_{1},0)$ and $(0,i_{2})$.\\
Let $M_{F}$ $\in \mathbf{M}_{2}(A/\chi _{\Phi })$ the  matrix of
the endomorphisme  Frobenius $F$ in this basis. Then
$M_{F}=~\left(
\begin{array}{ll}
a & b \\
a_{1} & b_{1}
\end{array}
\right)$, where $a,b,a_{1},b_{1}\in A/\chi _{\Phi }.$\\
 Although since :
Tr $M_{F}=a+b_{1}=c$ and $M_{F}(i_{1},0)=(i_{1},0)$ and $
M_{F}(0,i_{2})=(0,i_{2})$, we will have  $a.i_{1}\simeq i_{1}(\mod \chi_{\Phi })$ and then  $a-1$ is divisible by
$i_{1}$, of same for  $ b_{1}.i_{2}\simeq i_{2}(\mod \chi _{\Phi})$, that means that  $ b_{1}-1$ is divisible by $i_{2}$ and
then: $c-2=a-1+b_{1}-1$ is divisible by $i_{2}$ ( since we have
always  $i_{2}\mid i_{1}$).

Let \ $\rho $ be a prime ideal from $A$, different from the $A$
-characteristic $P$, we define the finite $A$-module $\Phi (\rho
)$ as been the $A$-module $(A/\rho )^{2}.$

The discriminant of the $A$-order: $A+g.O_{K(F)}$ is $\Delta .g^{2}$, where $%
\Delta $ is the discriminant of the characteristic polynomial
$P_{\Phi }(X)=X^{2}-cX + \mu P^{m}$. So each order is defined by
this discriminant and will be noted by $O($ disc$)$, see \cite
{Schoof}, and \cite{Tsfa-Vladut}. It is clear, by the
Propositions 2.1 that the inclusion $\Phi (\rho )\subset $
$L^{\Phi }$ implies that $\rho ^{2}\mid P_{\Phi }(1)$ and $\rho
\mid c-2$. We have :

\begin{proposition}
Let $\Phi $ be an ordinary Drinfeld $A$-module of rank 2, and let
$\rho $ an ideal from $A$ different from the $A$-characteristic $P$ of $L,$ such that $%
\rho ^{2}\mid P_{\Phi }(1)$ and $\rho \mid c-2$. Then $\Phi (\rho
)\subset L^{\Phi }$, if and only if, the $A$-order $\ O(\Delta
/\rho ^{2})\subset $ End$_{L}\Phi $.
\end{proposition}

To prove this proposition we need the following lemma :

\begin{lemma}
$\Phi (\rho )\subset L^{\Phi }$ is equivalent to
$\frac{F-1}{\varrho }\in $ End$_{L}\Phi $.
\end{lemma}

Proof : We know that $L^{\Phi }$ is satble by the isogeny $F$
  so $L^{\Phi }=$ Ker$(F-1)$, and by definition  $\Phi (\rho )$ $=$ Ker$(\rho )$ ( we
confuse by commodity the ideal $\rho $ with this generator in
$A$), and we know by \cite{Gosse}, Theorem 4.7.8, that for two
isogenies, let by example $F-1$ and $\rho $, we have
 Ker$(F-1)\subset $ \
Ker$(\rho )$, if and only if, there exists an element $g\in $ End$%
_{L}\Phi $ such that $F-1=g.\rho $ and then $\Phi (\rho )\subset
L^{\Phi },$ if and only if, $\frac{F-1}{\varrho }=g\in
$End$_{L}\Phi .$

We prove now the Proposition 2.3 :

Proof: Let $N(\frac{F-1}{\rho })$ the norm of the isogeny
$\frac{F-1}{\rho }$,  which is a principal ideal generated by
$\frac{P_{\Phi }(1)}{\rho ^{2}}$, and the trace $($Tr$)$ of this
isogeny is $\frac{c-2}{\rho }$ then we can calculate the
discriminant of the  $A$-module $ A[\frac{ F-1}{\rho }]$ by:

disc$A([\frac{F-1}{\rho }])=Tr(\frac{F-1}{\rho
})^{2}-4N(\frac{F-1}{\varrho } )=\frac{c^{2}-4\mu P^{m}}{\rho
^{2}}=\Delta /\rho ^{2}\Rightarrow$

$$O(\Delta /\rho ^{2})\subset  \ End_{L}\Phi.$$

We suppose now that : $O(\Delta /\rho ^{2})\subset $ End$_{L}\Phi
$ and we prove that $\Phi (\rho )\subset L^{\Phi }$. The Order
corresponding of the discriminant  $ \Delta /\rho ^{2}$ is
$A[\frac{F-1}{\rho }]$ this means that: $\frac{ F-1}{\varrho }\in
$ End$_{L}\Phi $ and so, by lemma 2.1 : $\Phi (\rho )\subset
L^{\Phi }$ .

\begin{corollary}
If $O(\Delta /\rho ^{2})\subset $ \ End$_{L}\Phi $, then $L^{\Phi
}$ is not cyclic.
\end{corollary}

Proof: We know that $\Phi (\rho )$ is not cyclic  (since it is a
$A$-module of rank  $2$), and then the necessary and sufficient
conditions need for non cyclicity of $A$-module \ $L^{\Phi }$ are
equivalent to the necessary and  sufficient conditions to have
$\Phi (\rho )\subset L^{\Phi }$.

We can so prove the following important theorem :

\begin{theorem}
Let $M=\frac{A}{I_{1}}\oplus \frac{A}{I_{2}}$, $I_{1}=(i_{1})$
And $\ \ I_{2}=(i_{2})$\, such that : $i_{2}\mid i_{1},$\
$i_{2}\mid (c-2)$. Then there exists an ordinary Drinfeld
$\mathit{A}$-module $\Phi $ over $L$ of rank $2$, such that:
$L^{\Phi }\simeq M$.
\end{theorem}

Proof : In fact, if we consider the Drinfeld $A$-module $\Phi$,
for which the characteristic of Euler-Poincare is giving by $\chi
_{\Phi }=I_{1}.I_{2}$ and his endomorphism ring is $O(\Delta
/i_{2}^{2})$ where $\Delta $ is always the discriminant of the
characteristic polynomial of the Frobenius $F$. We remind   that\
$\Phi (\rho )\subset L^{\Phi }$ for every  $\rho $ an ideal $A$,
different from $P$ and verify  $\rho ^{2}\mid P_{\Phi }(1)$ and
$\rho \mid (c-2),$ if and only if, the $A$-order $\ O(\Delta
/\rho ^{2})\subset $ \ End$_{L}\Phi $. Let now  $\rho =i_{2}$.
Since by construction the $A$-order $\ O(\Delta
/i_{2}^{2})\subset $ End$_{L}\Phi $ we have that $\Phi
(i_{2})\simeq (A/i_{2})^{2}\subset L^{\Phi }$. We know that
$L^{\Phi }$ is included or equal to $\Phi (\chi _{\Phi })$
$\simeq $ $\frac{A}{\chi _{\Phi }} \oplus \frac{A}{\chi _{\Phi
}}$, we have so  : $L^{\Phi }=\frac{ A}{I_{1}}\oplus
\frac{A}{I_{2}}$.

The above theorem can be proved by using the following conjecture:

\begin{conjecture}
Let $M \in \mathbf{M}_{2}( A / \chi_{\Phi })$, $\overline{P}= P($
mod $ \chi_{\Phi })$. We suppose : ( $det M$)=$\overline{P}^{m}$,
$Tr (M)=c$ and

$c$ $\nmid $ $P$.
There exists an ordinary Drinfeld $A$-module over a finite field $L$ of rank $%
2$, for which the Frobenius matrix associated, is $M_{F}$, and
such that:
$M_{F}= M \in \mathbf{M}_{2}(A/ \chi_{\Phi })$.
\end{conjecture}

We put the following matrix :
$M_{F}=\left( \begin{array}{cc} c-1       &  i_{1} \\i_{2}    &  -1\end{array} \right) \in \mathbf{M}_{2}(A/ \chi_{\Phi }).$

We can see that the three conditions of the conjecture are
realized then there exists an ordinary Drinfeld $A$-modules $\Phi
$ over $L$ of rank $2$, such that : $L^{\Phi }\simeq M$.

\selectlanguage{english}

\end{document}